\documentclass[12pt]{article}

\usepackage{amsmath}
\usepackage{amssymb}
\usepackage{amsfonts}
\usepackage{amsthm}
\usepackage{enumerate}

\newtheorem{theorem}{Theorem}[section]
\newtheorem{proposition}[theorem]{Proposition}
\newtheorem{lemma}[theorem]{Lemma}
\newtheorem{definition}[theorem]{Definition}

\newcommand{\rr}{{\mathbb R}}

\newcommand{\nn}{{\mathbb N}}

\newcommand{\calv}{{\cal V}}

\newcommand{\cals}{{\cal S}}

\linethickness{.4mm}

\begin{document}

\title{Incidence theorems for pseudoflats}

\author{Izabella {\L}aba and J\'ozsef Solymosi
\thanks{Both authors are supported in part by NSERC Discovery Grants.}
}

\maketitle 

\begin{abstract}
We prove Pach-Sharir type incidence theorems for a class of curves
in $\rr^n$ and surfaces in $\rr^3$, which we call {\em
pseudoflats}.  In particular, our results apply to a 
wide class of generic irreducible real algebraic sets
of bounded degree. 
\end{abstract}
 
\section{Introduction}\label{intro}

One of the most intriguing, and most actively studied, problems
in combinatorial geometry is finding upper bounds
on the number of point-curve and point-surface incidences.
The best known such result, which has since become an indispensable tool
with a wide variety of applications in discrete and combinatorial geometry,
is the Szemer\'edi-Trotter theorem \cite{SzT83} on point-line incidences
in the plane.  

There are now many extensions and generalizations of the Szemer\'edi-Trotter
theorem.  In one direction, Pach and Sharir \cite{PS98} have obtained
an analogous incidence bound for {\em pseudolines}, i.e. planar curves which
are uniquely determined by a certain fixed number $r$ of points (for lines,
we have $r=2$).  Another line of work concerns point-surface, especially
point-hyperplane, incidences in higher dimensions; see \cite{PS04} for
an excellent survey.  Recently, Elekes and T\'oth \cite{ET05} obtained
a sharp Szemer\'edi-Trotter type bound for point-hyperplane incidences in
$\rr^n$, assuming a certain non-degeneracy condition; this bound was
refined further by Solymosi and T\'oth \cite{ST05} under the additional
assumption that the point set in question is {\em homogeneous} (see
below).

The purpose of this paper is to propose a common generalization of the
results of \cite{PS98} and \cite{ET05}, \cite{ST05}, for homogeneous point
sets in $\rr^3$.  Specifically, we first obtain Pach-Sharir type
bounds for a class 
of curves (which we will also call pseudolines) in $\rr^n$, recovering 
a special case of the bound of \cite{PS98} for $n=2$.  We then use it to
prove the main result of our paper, namely an analogous incidence theorem
for a class of 2-dimensional surfaces in $\rr^3$
({\em pseudoplanes}).  
For 2-dimensional planes in $\rr^3$, our bound differs from that of
\cite{ET05}, \cite{ST05} only by the additional logarithmic factors; 
on the other hand, our non-degeneracy assumption is weaker than that
of \cite{ET05}, \cite{ST05}.

\begin{definition}\label{def-homogeneous}
(cf. \cite{SV04}, \cite{ST05}, \cite{Io1})
A finite point set $P\subset\rr^n$ is called {\em homogeneous} if 
$P$ lies in the interior of a $d$-dimensional cube $Q=[0,a]^n$ of volume
$\Theta(|P|)$\footnote{
Here and below, all constants depend only on $n$ and $r$.},
and if
any unit cube in $\rr^n$ contains at most $O(1)$ points of $P$.
\end{definition}

Fix $P$ as above.  We will say that a set $S\subset\rr^n$ (usually a curve
or surface) is {\em $m$-rich} if it contains at least $m$ points of $P$.
As in \cite{ET05}, \cite{ST05}, our results are stated in terms of a bound
on the number of $m$-rich pseudolines and pseudoplanes for a fixed point
set $P$ of cardinality $N$.  

Our first result concerns incidences for one-dimensional {\em pseudolines},
which we now define.

\begin{definition}\label{def-pseudolines}
(cf. \cite{PS98})
Let $\calv$ be a family of subsets in $\rr^n$.  We say that $\calv$ is a {\em type
$r$ family of pseudolines} if the following two conditions are satisfied:

\smallskip
(i) {\em (rectifiability)} Let $t\in\nn$.  If the enclosing cube $Q$ is subdivided 
into $t^n$ congruent and disjoint (except for boundary) subcubes, then each
$V\in\calv$ has nonempty intersection with at most $O(t)$
 subcubes.

\smallskip
(ii) {\em (type $r$)} For any distinct $r$ points in $\rr^n$, there is at most
one $V\in\calv$ which contains them all.
\end{definition}

In Section \ref{algebraic}, we verify that the
conditions of Definition \ref{def-pseudolines} hold if $\calv$ is a family 
of irreducible one-dimensional algebraic varieties defined by polynomial 
equations of degree at most $d$ (with an explicit value of $r$, depending
on $n$ and $d$).
It is useful to think of the elements of $\calv$ as curves, since
this description applies to the main cases of interest (such as 
the algebraic varieties just mentioned).  Note that
the rectifiability condition (ii) implies that each $V\in\calv$ has Hausdorff
dimension at most 1.  However, no
continuity or smoothness assumptions are actually required.

\begin{theorem}\label{thm-1dim}
Let ${\cal V}$ be a type $r$ family of pseudolines in 
$\rr^n$, and let $P$ be a homogeneous set of $N$ points.  
Then:

(i) if $k\geq CN^{1/n}$ for large enough $C$, then there are no 
$k$-rich pseudolines in $\calv$;

(ii) if $k\leq CN^{1/n}$, then
the number of $k$-rich pseudolines in $\calv$ is bounded by
$O(N^r/k^{n(r-1)+1}).$
\end{theorem}

For $n=2$, we recover a special case (for homogeneous point sets) 
of the Pach-Sharir theorem on incidences
for pseudolines \cite{PS98}.  For $n\geq 3$, Theorem \ref{thm-1dim} extends a 
result of Solymosi and Vu \cite{SV04} on incidences for lines in $\rr^n$.  
Our proof is in fact very similar to that of \cite{SV04}.

Our main result concerns 2-dimensional surfaces in $\rr^3$.  (Again, the
rectifiability assumption (i) implies Hausdorff dimension at most 2, but
otherwise our ``surfaces" could be quite arbitrary.)

\begin{definition}\label{def-pseudoplanes}
Let $\cals$ be a family of subsets of $\rr^3$.  We say that $\calv$ is a {\em type
$r$ family of pseudoplanes} if the following holds:

\smallskip
(i) {\em (rectifiability)} Let $t\in\nn$.  If the enclosing cube $Q$ is subdivided 
into $t^3$ congruent and disjoint (except for boundary) subcubes, then each
$S\in\cals$ has nonempty intersection with at most $O(t^2)$ subcubes.

\smallskip
(ii) {\em (pairwise intersections are pseudolines)} 
Let
$$\calv=\{S\cap S':\ S,S'\in\cals,\ S\neq S'\}.$$
Then $\calv$ is a type $r$ family of pseudolines.

\end{definition}

\begin{theorem}\label{surf-thm}
Let $P$ be a homogeneous point set in a cube $Q$, $|P|=N$.
Let $\cals$ be a type $r$ family of pseudoplanes in $\rr^3$.
Assume that each $S\in\cals$ contains a {\em defining $r+1$-tuple},
i.e.~$r+1$ distinct points $x_1,\dots,x_{r+1}$ such that no pseudoplane
in $\cals$, other than $S$ itself, contains them all.  
Then:

\smallskip
(i) if $k\geq CN^{2/3}$ for large enough $C$, then there are no 
$k$-rich pseudoplanes in $\cals$;

\smallskip
(ii) if $k\leq CN^{2/3}$, then
the number of $k$-rich pseudoplanes in $\cals$ is bounded by
$O(N^{r+1}(\log N\log k)^{3r/2+2}/k^{3r/2+1}).$
\end{theorem}

{\bf Remark.} It is well known (cf. \cite{ET05}, Section 2) that
bounds on the number of $k$-rich curves or surfaces are equivalent
to the more standard formulation in terms of a bound on the total
number of incidences between the point set and the objects in
question.  Essentially (i.e. modulo the endpoints), we obtain
bounds $O(M^{\frac{n(r-1)}{n(r-1)+1}} N^{\frac{r}{n(r-1)+1}})$
on the number of incidences between $M$ type $r$ pseudolines and
$N$ well-distributed points in $\rr^n$, and 
$O(M^{\frac{3r}{3r+2}} N^{\frac{2r+2}{3r+2}})$
on the number of incidences between $M$ type $r$ pseudoplanes and
$N$ well-distributed points in $\rr^3$.

\medskip

We briefly explain how our results fit in with the existing literature on 
point-plane incidences.  It is clear that if $m$ points of $P$ lie on one line,
then any plane containing this line will be $m$-rich, and there is
no bound on the number of such planes.  Therefore any non-trivial
Szemer\'edi-Trotter type results for higher-dimensional surfaces must
make some non-degeneracy assumption on the plane-point configuration.
The special cases considered in the literature include configurations
where there are no three collinear points \cite{EGS}, the incidence
graph does not contain a $K_{r,r}$ \cite{BK}, and when all hyperplanes
are spanned by the point set \cite{AA} (see also \cite{EdS}).  

Our result is closest to those of \cite{ET05} and \cite{ST05}. Elekes
and T\'oth \cite{ET05}
give a sharp bound $O(\frac{N^n}{k^{n+1}}+\frac{N^{n-1}}{k^{n-1}})$
on the number of $k$-rich $n-1$-dimensional hyperplanes in $\rr^n$,
with respect to a point set $P$
of cardinality $N$, provided that the hyperplanes are {\em not-too-%
degenerate} in the following sense: there is an $\alpha<1$ such that 
for each hyperplane $S$, no more than $\alpha|P\cap S|$ points of $P
\cap S$ lie in a lower-dimensional flat.  This estimate was strengthened
by Solymosi and T\'oth \cite{ST05} to $O(\frac{N^n}{k^{n+1}})$, under the same
non-degeneracy assumption and with the additional condition that the
point set $P$ is homogeneous; the homogeneity assumption is necessary
here, as Elekes and T\'oth observe that for certain ranges of $N,k$
the second term in the bound
of \cite{ET05} is dominant and sharp.

The bound in Theorem \ref{surf-thm} with 
$r=2$ matches that of \cite{ST05} for $n=3$, modulo the extra logarithmic
factors.  On the other hand, the nondegeneracy condition of \cite{ET05},
\cite{ST05} is stronger than that of Theorem \ref{surf-thm}, where 
it suffices for each plane to contain just one defining triple (i.e. 
non-trivial triangle) of points in $P$.  We do not know whether the 
logarithmic factors can be dropped without additional assumptions (such
as those in \cite{ET05}, \cite{ST05}) on the distribution of {\em all}
points in each plane.


If we only assume that each $n-1$-dimensional hyperplane contains $n$
affinely independent points of $P$, but no further conditions are imposed on
the hyperplane-point configuration, the bounds of \cite{ET05}, \cite{ST05}
are known to be false.  In this case, an optimal incidence bound is due to 
Agarval and Aronov \cite{AA}, namely there are at most $O(\frac{N^n}{k^3}
+\frac{N^{n-1}}{k})$ $k$-rich hyperplanes spanned by $P$.   Examples 
of \cite{Ed}, \cite{EdH} show that this estimate is sharp.  This also shows
that our theorem would already fail for $r=2$ without the homogeneity
assumption on $P$.  

Incidence bounds for curves and surfaces, other than lines and planes,
in dimensions 3 and higher, are in general not well understood.
Non-trivial bounds have been obtained only in certain special cases,
e.g. circles, spheres, algebraic curves of bounded degree and
lying in distinct planes.  See \cite{PS04} for an overview.

A very rough outline of the proof of Theorem \ref{surf-thm}
is as follows.  Consider a sequence of nested cell decompositions of
the enclosing cube $Q$.  In a coarse decomposition (e.g. when the
entire $Q$ is a single cell), there are many same-cell defining
$r+1$-tuples of points of $P$.  In a very fine
decomposition, e.g. if each cell contains only one point of $P$,
there are no such $r+1$-tuples.  We find an intermediate scale at
which the transition takes place for most points and surfaces.
At that scale, if each cell contains few points of $P\cap S$ for most
$S$, we simply double-count the number of defining same-cell
$r+1$-tuples.  If on the other hand each cell contains many points
of $P\cap S$, then these points must in fact live on the intersection 
pseudolines in $\calv$.  Our bound is now obtained by applying 
Theorem \ref{thm-1dim} to $\calv$ and to the points of $P$
in each cell.

Our argument does not seem to extend easily to hypersurfaces in higher
dimensions, and in fact it is not clear how one should define
higher-dimensional pseudoflats.  For example, our proof relies 
heavily on the assumption that pairwise intersections of surfaces 
are pseudolines defined uniquely by $r$ points; this condition
fails for generic higher-dimensional hypersurfaces where pairwise
intersections may have dimension 2 or more.  There seems to be no
easy way to circumvent this by considering multiple intersections.


The plan of the paper is as follows.  In Section \ref{algebraic}, we discuss
the applicability of Theorems \ref{thm-1dim} and \ref{surf-thm} to 
algebraic sets, and give a few examples to illustrate this.  We then
prove Theorems \ref{thm-1dim} and \ref{surf-thm} in Sections
\ref{1-dim} and \ref{sec-surfaces}, respectively.

We use Roman letters to denote Cartesian coordinates of points in $\rr^n$, eg.
$x=({\rm x}_1,\dots,{\rm x}_n)\in\rr^n$.


\section{Algebraic varieties as pseudoflats}
\label{algebraic}

In this section we state conditions on families of algebraic varieties
under which the assumptions of Definitions \ref{def-pseudolines} and
\ref{def-pseudoplanes} are satisfied.

We briefly recall a few basic definitions from algebraic geometry,
restricting our attention to {\em real} algebraic sets.
The reader is cautioned that the terminology and features of real
algebraic geometry are sometimes quite different from the complex
case; see e.g. \cite{BCR} for more details.

For the purposes of this paper, 
a {\em real algebraic set} or {\em real algebraic variety} 
(we will usually omit the qualifier ``real" in what follows)
is the zero set in $\rr^n$ of a finite family of polynomials $F_1(x),\dots,
F_s(x)$ with real coefficients\footnote{
In real algebraic geometry, one also considers algebraic sets over
{\em real closed fields}.  Furthermore, it is useful to distinguish
between an algebraic set and an algebraic variety, the
latter being an algebraic set equipped with a sheaf of regular
functions.  However, we do not need to make this distinction here.}.
We will say that an algebraic set $S$
is {\em reducible} if there are two algebraic sets $S',S''$, neither
equal to $S$, such that $S=S'\cup S''$; otherwise, we will say that $S$ is 
{\em irreducible}.

There are several equivalent ways of defining the {\em dimension} of
an algebraic set.  The easiest one for us to use is the
following: the dimension $k$ of an algebraic set $S$ is the length of
the longest chain of irreducible varieties $S_j$ such that
$$\emptyset\neq S_1\subsetneq S_2\subsetneq \dots \subsetneq S_k\subset S.$$
In particular, an irreducible variety does not contain any proper
subvariety of the same dimension.  

If an algebraic set $S$ is a $C^\infty$ $k$-dimensional submanifold of 
$\rr^n$, then its algebraic dimension is $k$.  Note, however, that there
are irreducible real algebraic sets which consist of several components
of different topological dimensions (see e.g. the examples in 
\cite{BCR}, pp. 60--61).  In such cases, the algebraic dimension
of the set will be the largest of the dimensions of its components.

\begin{proposition}\label{alg-dim1}
Let ${\cal V}$ be a family of irreducible one-dimensional varieties in $\rr^n$,
defined by a polynomial equations of degree at most $d$.
Then ${\cal V}$ is a type $r$ family of pseudolines, with $r=d^2+1$ if
$d=2$, and with $r= d(2d-1)^{n-1}+1$ if $n\geq 3$.
\end{proposition}

\begin{proposition}\label{alg-dim2}
Let $\cals$ be a family of 2-dimensional algebraic varieties in $\rr^3$,
each given by a polynomial equation of degree no more than $d$.
Assume that the intersection $V=S\cap S'$ of any two varieties 
$S,S'\in \cals$, $S\neq S'$,  is an irreducible one-dimensional variety.
Then $\cals$ is a type $r$ family of pseudoplanes, with $r=
d(2d-1)^{2}+1$.
\end{proposition}

The proofs of both propositions will rely on the following
result from real algebraic geometry
\cite{BPR96}, \cite{BPR05}.   
Let $V\subset \rr^n$ be a $k$-dimensional variety 
defined by polynomials of degree at most $d$.  Let also $P_1,
\dots,P_s$ be polynomials in $n$ variables of degree at most $d$.
A {\em sign condition} for the set ${\cal P}=\{P_1,\dots,P_s\}$
is a vector $\sigma\in\{-1,0,1\}^s$.  
We write
$$
\sigma_{{\cal P},V}=\{x:\ x\in V,({\rm sign}(P_1(x)),\dots,
{\rm sign}(P_n(x)))=\sigma\},
$$ 
and call its non-empty semi-algebraically connected components
{\em cells} of the sign condition $\sigma$ for ${\cal P}$ over 
$V$. Let $|\sigma_{{\cal P},V}|$ be the number of such cells, 
then
$$
C({\cal P},V)=\sum_\sigma |\sigma_{{\cal P},V}|
$$
is the number of all cells defined by all possible sign conditions.
Let $f(d,n,k,s)$ be the maximum of $C({\cal P},V)$ over all
varieties $V\subset\rr^d$ and sets of polynomials ${\cal P}$ as
described above.
Then the main result of \cite{BPR96} (see also \cite{BPR05}) is
that
\begin{equation}\label{bpr}
f(d,n,k,s)={s\choose k}(O(d))^n.
\end{equation}

\medskip

{\it Proof of Proposition \ref{alg-dim1}.}
If $\calv$ is a family of one-dimensional irreducible varieties in $\rr^2$,
each defined by a polynomial equation of degree at most $d$, it follows
from Bezout's theorem that any two distinct varieties in $\calv$ intersect
in no more than $d^2$ points, hence $\calv$ is type $r$ for $r=d^2+1$.  

Suppose now that $\calv$ is a family of one-dimensional irreducible varieties
in $\rr^n$, $n\geq 3$, each one defined by a system of polynomial equations
in $n$ variables of degree at most $d$.  The classic results on the sum of
Betti numbers of algebraic sets \cite{OP49}, \cite{M64}, \cite{Thom65} imply
in particular that the intersection of two such distinct varieties has no more than 
$d(2d-1)^{n-1}$ connected components; since the varieties are irreducible,
the intersection is a variety of dimension 0, hence each connected 
component is a single point.  Thus we may take
$r=d(2d-1)^{n-1}+1$.

It remains to verify rectifiability.
Let $t\in \nn$.
We subdivide the enclosing cube $Q$ into $t^n$ congruent open subcubes
$Q_j=\{(j_i-1)a/t<{\rm x}_i<j_ia/t,\ i=1,\dots,n\}$, indexed by
$j=(j_1,\dots,j_n)\in\{1,\dots,t\}^n$.
Fix a $V\in\calv$, and consider the sets 
$V_j=V\cap Q_j$. We have $P\cap V\subset\bigcup_j(P\cap
V_j)$.  Each nonempty $V_j$ contains a cell of $V$
associated with a suitable sign condition for the system of
polynomials $P_{i,s}(x)= {\rm x}_i-sa/t$, $i=1,\dots,n$,
$s=1,\dots,t$.  By (\ref{bpr}), 
the number of such cells is bounded by $ntO(d)^n$, as required.
\qed

\medskip

{\it Proof of Proposition \ref{alg-dim2}.}
Let $\calv=\{S\cap S':\ S,S'\in\cals,\ S\neq S'\}$, then
$\calv$ is a type $r$ family of pseudolines by Proposition \ref{alg-dim1}.
The proof of rectifiability is the same as in the proof of
Proposition \ref{alg-dim1}, except that this time (\ref{bpr}) yields
the bound ${nt\choose 2}O(d)^n$ on the number of non-empty cells. 
\qed

\medskip


\section{Proof of Theorem \ref{thm-1dim}}
\label{1-dim}

We assume that the enclosing cube for $P$ is $Q=[0,a]^n$
for some positive
integer $a$, and that all points in $P$ have irrational coordinates.
We always let $N$ be sufficiently large.  Without loss of generality,
we assume that all $V\in\calv$ are $k$-rich.

We first prove (i).  Let $t$ be an integer to be fixed later.  
We subdivide the enclosing cube $Q$ into $t^n$ congruent open subcubes
$Q_j=\{(j_i-1)a/t<{\rm x}_i<j_ia/t,\ i=1,\dots,n\}$, indexed by
$j=(j_1,\dots,j_n)\in\{1,\dots,t\}^n$.
By the assumption from the last paragraph, no points in $P$ lie on the boundary
of any $Q_j$.  
For each $V\in\calv$, we let
$V_j=V\cap Q_j$. We have $P\cap V\subset\bigcup_j(P\cap
V_j)$.  By the rectifiability assumption, the number of non-empty $V_j$'s
is at most $O(t)$.

We now choose $t$ so that $t^n=\Theta(N)$.  
Since $P$ is homogeneous, each ${\cal Q}_j$ contains no more than
$O(1)$ points of $P$.  Hence the cardinality of $P\cap V$ is bounded
by $O(t)=O(N^{1/n})$, as claimed.

It remains to prove (ii).  We divide
$Q$ into $t^n$ subcubes $Q_j$ as in the proof of (i), except
that this time we will choose 
\begin{equation}\label{e-t}
t=\Theta(k),
\end{equation}
with the implicit constants small enough (depending on $n,r$).
In particular, by (i) we may assume that $t^n\leq N$, since otherwise
there is nothing left to prove.

A $r$-tuple of distinct points $x_1,\dots,x_r\in P$ is {\em good} if all
$x_i$ belong to the same subcube $Q_j$.  
We count the number $M$ of good $r$-tuples in two ways.  On one hand, 
since $P$ is homogeneous, each $Q_j$ contains no more than
$O(N/t^n)$ points of $P$.  Thus
$$M= O((N/t^n)^r\cdot t^n)=O(N^r/t^{n(r-1)}).$$

On the other hand, let $V_j=V\cap Q_j$ as in the proof of (i).
Each $r$-tuple of points in $P\cap V$ contained in
one $P\cap V_j$ is good.  By rectifiability, the
number $K$ of distinct and non-empty $V_j$'s is bounded by
$k/r$ (provided that the constants in (\ref{e-t}) were chosen
appropriately small).  Thus the number of good $r$-tuples in $P\cap V$ is 
$\Omega((k/K)^r\cdot K)=\Omega(k^r/K^{r-1})=\Omega(k)$.

Summing over all $V\in \calv$ and remembering that any $r$-tuple can
belong to only one $V$ (since $\calv$ is type $r$), 
we see that
$$M=\Omega(|\calv|k).$$
Comparing the upper and lower bounds for $M$, and using (\ref{e-t}), 
we see that
$$|\calv|=O(\frac{N^r}{kt^{n(r-1)}})
=O(\frac{N^r}{k^{n(r-1)+1}})$$
as claimed.

\section{Proof of Theorem \ref{surf-thm}}
\label{sec-surfaces}

The proof of (i) is identical to that of Theorem \ref{thm-1dim}(i),
except that the 2-dimensional rectifiability condition
yields the exponent $2/n$ as indicated in the theorem.  We omit the details.

We now prove (ii).  Let $|\cals|=X$.  
We assume that $Q=[0,a]^3$ for some positive
integer $a$, and that all points in $P$ have irrational coordinates.

For $i=0,1,2,\dots,I$, we define the $i$-th cutting of $Q$ to be the
subdivision of $Q$ into $2^{3i}$ congruent open subcubes
$Q_{i,j}$ of sidelength $a/2^i$.  We let $I=\Theta(\log N)$
so that each subcube in the $I$-th cutting contains at most 1 point of $P$.
Note that no points in $P$ lie on the boundary of any $Q_{i,j}$.  

For each $i$ and each $S\in\cals$, the $i$-th cutting divides $S$ into
subsets $S_{i,j}=S\cap Q_{i,j}$.  By the rectifiability assumption,
we have
\begin{equation}\label{cell-bound}
|\{j:\ S_{i,j}\neq\emptyset\}|=O(2^{2i}),
\end{equation}
with constants uniform in $i$.

Let $S\in\cals$.  We will say that a $r+1$-tuple of points $x_1,\dots,
x_{r+1}$ is {\em defining for $S$ at level $i$} if $x_1,\dots,x_{r+1}$
are distinct points in $P\cap S$ which all belong to the same subcube
of the $i$-th cutting, and if moreover
there is no other surface $S'\in\cals$, $S'\neq S$,
such that $x_1,\dots,x_{r+1}\in S'$.  Thus a defining $r+1$-tuple at
level 0 is simply a defining $r+1$-tuple for $S$ as in the statement 
of the theorem.

We define the {\em index} $i(x,S)$ of a pair $(x,S)$, where
$S\in \cals$, $x\in P\cap S$, to be the least value of
$i$ such that $x$ does not belong to a
defining $r+1$-tuple for $S$ at level $i$.  

\begin{lemma}\label{index-lemma}
For all $S\in \cals$, $x\in P\cap S$, we have $1\leq i(x,S)\leq I$.
\end{lemma}

{\it Proof.}
Clearly there are no defining $r+1$-tuples at 
level $I$, hence $i(x,S)\leq I$.  
It remains to prove that $i(x,S)\geq 1$ for all $x\in P\cap S$, $S\in\cals$.  
Indeed, fix $S\in\cals$ and $x_0\in P\cap S$.  We need to prove that 
$x_0$ belongs to a defining $r+1$-tuple for $S$.  By the non-degeneracy 
assumption, $S$ contains a defining $r+1$-tuple $T=\{x_1,\dots,x_{r+1}\}$. 
If $x_0\in T$, we are done.  Otherwise, let $T_j=(T\setminus \{x_j\})\cup
\{x_0\}$, and suppose that $T_1,\dots,T_{r+1}$ are all non-defining.  
This means that for each $j=1,\dots,r+1$ there is a $S_j\in\cals$, 
$S_j\neq S$, such that $T_j\subset V_j:=S\cap S_j$.  In particular,
if $j\geq 2$, then $V_1$ and $V_j$ share the $r$ points $x_m$, 
$m\in\{0,1,\dots,r+1\}\setminus\{1,j\}$.  By the $r$-type assumption,
$V_1=V_j$, $j=2,\dots,r+1$.  But then all the points $x_1,\dots,x_{r+1}$
belong to $V_1=S\cap S_1$, hence $T$ is not defining, contradicting our
assumption.  It follows that at least one of $T_1,\dots,T_{r+1}$ is
a defining $r+1$-tuple for $S$ containing $x_0$, as required.
\qed

\medskip

For each $S$, we choose $i(S)$ to be the least value of
$i$ such that
$$
|\{x:\ x\in P\cap S,\ i(x,S)=i(S)\}|\geq k/2I.
$$
We then choose an ${\bf i}\in\{0,1,\dots,I\}$ and a subset
$\cals_1\subset \cals$ such that
$$
|\cals_1|\geq|\cals|/2I,\ i(S)={\bf i}
\hbox{ for all } S\in\cals_1.
$$

\medskip

{\bf Case 1:} $k\leq C_0 2^{2{\bf i}}\log N\log k$.  
We count the number
$M$ of all defining $r+1$-tuples for all $S\in\cals_1$ at level
${\bf i}-1$.  Each $S\in\cals_1$ contains at least $k/2I$ points
$x\in P$ with index $i(x,S)={\bf i}$.  Each such point must belong to a 
defining $r+1$-tuple for $S$ at level ${\bf i}-1$, and each $r+1$-tuple
can be defining for only one $S$.  Thus 
$$
M\geq \frac{X}{2I}\cdot\frac{k}{2I}\cdot\frac{1}{r+1}.
$$
On the other hand, $M$ is trivially bounded from above by 
the total number of the $r+1$-tuples that belong to the same 
cube of the ${\bf i}-1$-th cutting, 
$$
M=O\Big((\frac{N}{2^{3{\bf i}-3}})^{r+1}\cdot 2^{3{\bf i}-3}\Big)
=O(N^{r+1}/2^{3{\bf i}r}).
$$
Comparing the upper and lower bounds, and using the assumption
on $k$ for Case 1, we get
$$
X=O\Big(\frac{N^{r+1}I^2}{k\cdot 2^{3{\bf i}r}}\Big)
=O\Big(\frac{N^{r+1}(\log N\log k)^{3r/2}(\log N)^2}{k^{3r/2+1}}\Big)
$$
as required.

\medskip
{\bf Case 2:} $k\geq C_0 2^{2{\bf i}}\log N\log k $.
In this case, points of $P$ tend to be aligned along the one-dimensional
intersection curves; we will therefore use our one-dimensional 
incidence bound.  We first do some pigeonholing to fix the values of 
certain parameters.  
For each $S\in\cals_1$, we let
$$
P(S)=\{x\in P\cap S:\ i(x,S)={\bf i}\},
$$
then $|P(S)|\geq k/2I$.  We then choose a subset $P_0(S)
\subset P(S)$ such that $|P_0(S)|\in[\frac{k}{2I},\frac{k}{2I}+1)$.
Let $L$ be an integer such that $2^L\leq k<2^{L+1}$ (hence $L=\Theta(\log k)$).
Note that for each $j$,
$$ |P_0(S)\cap S\cap {\cal Q}_{{\bf i},j}|
\leq |P_0(S)|\leq \frac{k}{2I}+1\leq k+1\leq 2^{L+1}.$$
Thus if we let
$$
m(l,S)=|\{j:\ |P_0(S)\cap S\cap Q_{{\bf i},j}|\in
[2^l,2^{l+1}]\}|,\ l=0,1,\dots,L,
$$
then for each $S\in\cals_1$,
$$
\sum_{l=0}^{L} m(l,S)\cdot 2^{l+1}\geq |P_0(S)|\geq k/2I,
$$
hence we may choose $l(S)$ such that 
$$
m(l(S),S)\cdot 2^{l(S)+1}\geq k(4IL)^{-1}.
$$
Pigeonholing again, we find
a value of ${\bf l}\in\{0,\dots,L\}$ and a set $\cals_2\subset
\cals_1$ such that
$$
|\cals_2|\geq|\cals_1|/{2L},\ 
l(S)={\bf l}\hbox{ for all }S\in\cals_2.
$$

Let $S\in\cals_2$, and let $S_j=S\cap Q_{{\bf i},j}$.  Relabelling
the subcubes if necessary, we may assume that
$$
|P_0(S)\cap S_j|\geq 2^{\bf l},\ j=1,\dots, m,
$$
where 
\begin{equation}\label{m-vs-k}
m\cdot 2^{{\bf l}+1}\geq k(4IL)^{-1}.  
\end{equation}
By rectifiability, we have $m\leq C\cdot 2^{2{\bf i}}$.  Thus it follows that
\begin{equation}\label{enough-points}
2^{{\bf l}}\geq k(8ILm)^{-1}\geq r+1,
\end{equation}
provided that the constant $C_0$ in the assumption of Case 2 was
chosen large enough.

We now claim that for each $j=1,\dots,m$, there is a unique 
$V_j\in\calv$ such that
$P_0(S)\cap S_j\subset V_j.$  Indeed, let $x_1,\dots,x_r,x_{r+1}\in
P_0(S)\cap S_j$.  By the definition of $P_0(S)$ and ${\bf i}$,
$x_1,\dots,x_{r+1}$ is not a defining $r+1$-tuple for $S$, hence
there is a $S'\in \cals$, $S'\neq S$ such that $x_1,\dots,x_{r+1}
\in S'$.  Thus $x_1,\dots,x_{r_1}\in V_j:=S\cap S'$.  Let now
$x\in P_0(S)\cap S_j$, $x\neq x_1,\dots,x_{r+1}$.  Then $x_1,
\dots,x_r,x$ is another non-defining $r$-tuple, hence 
$x_1,\dots,x_{r},x\in V'_j:=S\cap S''$ for some $S''\in\cals$.
But then $V'_j$ intersects $V_j$ in $r$ distinct points $x_1,
\dots,x_r$, hence $V'_j=V_j$ since $\calv$ is type $r$.
It follows that $x\in V_j$ for all $x\in P_0(S)\cap S_j$, as
claimed.

Thus for each $S\in \cals_2$, there are at least $m$ subcubes 
$Q_{{\bf i},j}$ with the following property: 
there is a $V_j=V_j(S)\in\calv$ which contains at least $2^{2{\bf l}}$
points of $S\cap Q_{{\bf i},j}$ with index $i(x,S)={\bf i}$.
Moreover, we have (\ref{m-vs-k}) and (\ref{enough-points}).

Now for the main argument.  We count the number $M'$ of ``admissible"
triples $(S,V,j)$ such that:
\begin{itemize}
\item $S\in\cals$, $V\in\calv$, $V=S\cap S'$ for some $S'\in\cals$;
\item $|V\cap P\cap Q_{{\bf i},j}|\geq 2^{\bf l}.$
\end{itemize}

\noindent
{\bf Lower bound:} From the above construction, for each $S\in\cals_2$
there are at least $m$ values of $j$ such that 
$$|V\cap P\cap Q_{{\bf i},j}|\geq 
|P_0(S)\cap S_j|\geq 2^{\bf l}$$
for some $V\in\calv$ (depending on $j$).  Hence
$$
M'\geq |\cals_2|\cdot m
=\Omega(X\cdot k/IL).
$$

\noindent
{\bf Upper bound:}  There are three ingredients.

\begin{itemize}
\item First, there are at most $2^{3{\bf i}}$ values of $j$.  

\item For each $j$, we estimate the number of eligible $V$'s by 
applying Theorem \ref{thm-1dim} to $\calv$ and to point
sets $P_j:=P\cap Q_{{\bf i},j}$, homogeneous in
$Q_{{\bf i},j}$ and of cardinality $\Theta(N/2^{3{\bf i}})$.
Thus the number of $V\in\calv$ containing at least
$2^{\bf l}$ points of $P_j$ is bounded by
$$
O\Big(\frac{(N/2^{3{\bf i}})^r}{(2^{\bf l})^{3r-2}}\Big).
$$

\item Finally, we claim that for each such fixed $j$ and $V$, there
are at most $O(N/2^{3{\bf i}})$ surfaces $S\in\cals$ such that
$(S,V,j)$ is admissible.  Indeed, define
the {\em parent cube} of $Q_{{\bf i},j}$ to be the unique 
cube in the ${\bf i}-1$-th cutting which contains 
$Q_{{\bf i},j}$.  Suppose that $(S,V,j)$ is admissible.  
Then $V\cap Q_{{\bf i},j}$ contains at least $2^{\bf l}$
points $x\in P$ with $i(x,S)={\bf i}$.  Fix such an $x$, then by the
definition of index, $x$ belongs to a defining
$r+1$-tuple for $S$ at level ${\bf i}-1$, i.e. contained in the
parent cube.  This $r+1$-tuple must
contain at least one point, say $x_0$, which is not in $V$.
It remains to prove that $V$ and $x_0$ define $S$ uniquely; this
implies the claim, since the parent cube contains
at most $O(N/2^{3{\bf i}})$ points of $P$.

By (\ref{enough-points}), there are at least $r+1$ distinct points $x_1,
\dots,x_{r+1}$ in $V\cap P_j$.  It suffices to prove that $x_0,
x_1,\dots,x_r$ is a defining $r+1$-tuple for $S$ at level ${\bf i}
-1$.  Indeed, suppose to the contrary that there is a $S''\in\cals$,
$S''\neq S$, such that $x_0,x_1,\dots,x_r\in V':=S\cap S''$.  
But then $x_1,\dots,x_r\in V\cap V'$.  Since $\calv$ is type $r$, it
follows that $V'=V$, and in particular that $x_0\in V$, contrary to
our choice of $x_0$.  
\end{itemize}

Combining the three estimates, we obtain the upper bound
$$
M'= O\Big(\frac{(N/2^{3{\bf i}})^r}{(2^{\bf l})^{3r-2}}
\cdot\frac{N}{2^{3{\bf i}}}\cdot 2^{3{\bf i}}\Big)
= O\Big(\frac{N^{r+1}2^{-3{\bf i}r}}{(2^{\bf l})^{3r-2}}\Big).
$$

\noindent
{\bf Conclusion:}
Comparing the upper and lower bounds on $M'$, we get that
$$
X
= O\Big(\frac{N^{r+1}2^{-3{\bf i}r}\log N\log k}{(2^{\bf l})^{3r-2}m}\Big).
$$
By (\ref{m-vs-k}), we have $2^{\bf l}\geq k(m\log N\log k)^{-1}$.
Hence
$$
X
=O\Big(\frac{N^{r+1}2^{-3{\bf i}r}\log N\log k}{(k/m\log N\log k)^{3r-2}m} \Big)$$
$$=O\Big(\frac{N^{r+1}2^{-3{\bf i}r}m^{3r-3}}{k^{3r-2}} (\log N\log k)^{3r-1}\Big).$$
By rectifiability, we have $m=O(2^{2{\bf i}})$, so that
$$
X
=O\Big(\frac{N^{r+1}2^{(3r-6){\bf i}}}{k^{3r-2}} (\log N\log k)^{3r-1}\Big).$$
Finally, the assumption of Case 2 is that
$2^{2{\bf i}}=O(k/\log N\log k)$.  Thus
$$X
=O\Big( C\frac{N^{r+1}(k/\log N\log k)^{3r/2-3}}{k^{3r-2}} (\log N\log k)^{3r-1}\Big)$$
$$=O\Big(\frac{N^{r+1}}{k^{3r/2+1}} (\log N\log k)^{3r/2+2}\Big).$$
This completes the proof of the theorem.



\noindent{\sc Department of Mathematics, University of British Columbia, Vancouver, 
B.C. V6T 1Z2, Canada}

\noindent{\it ilaba@math.ubc.ca, solymosi@math.ubc.ca}

\end{document}